\documentclass[10pt,leqno]{article}
\usepackage{amssymb,amsbsy,amsmath,amsfonts,amssymb,amscd, mathrsfs}
\usepackage{graphicx} 

\usepackage[all]{xy}
\usepackage{color} 

\usepackage[english]{babel}
\usepackage[T1]{fontenc}
\usepackage{indentfirst}

\makeatletter
\@addtoreset{equation}{section}
\makeatother

\newtheorem{statement}{}[section]
\newtheorem{theorem}[statement]{Theorem}
\newtheorem{lemma}[statement]{Lemma} 

\newtheorem{proposition}[statement]{Proposition}

\newtheorem{corollary}[statement]{Corollary}

\newcommand\C{\mathbb C}

\newcommand\R{\mathbb R}
\newcommand\T{\mathbb T}
\newcommand\D{\mathbb D}
\newcommand\Z{\mathbb Z}
\renewcommand\H{\mathbb H}
\newcommand\e{{\rm e}}

\newcommand\eps{\varepsilon}
\newcommand\ind{{\rm 1\kern-.30em I}}
\newcommand\qed{\hfill $\square$}
\renewcommand \Re{{\mathfrak R}{\rm e}\,}
\renewcommand \Im{{\mathfrak I}{\rm m}\,}

\let\phi=\varphi

\title{\bf Compact composition operators on the Dirichlet space and capacity of sets of contact points}
\author{\it Pascal Lef\`evre, Daniel Li, Herv\'e Queff\'elec, \\ 
\it Luis Rodr{\'\i}guez-Piazza\footnote{Supported by a Spanish research project MTM 2009-08934.}}

\date{\footnotesize \today}

\begin{document}

\maketitle

\noindent{\bf Abstract.} \emph{We prove that for every compact set $K \subseteq \partial{\D}$ of logarithmic capacity ${\rm Cap} \, K = 0$, there exists a 
Schur function $\phi$ both in the disk algebra $A (\D)$ and in the Dirichlet space ${\cal D}_\ast$ such that the composition operator $C_\phi$ is in all Schatten  
classes $S_p ({\cal D}_\ast)$, $p > 0$, and for which $K = \{\e^{i t} \, ; \ |\phi (\e^{it})| = 1 \} = \{\e^{i t} \, ; \ \phi (\e^{it}) = 1 \}$. We show that  for 
every bounded composition operator $C_\phi$ on  ${\cal D}_\ast$ and every $\xi \in \partial \D$, the logarithmic capacity of 
$\{ \e^{it} \, ; \ \phi (\e^{it}) = \xi \}$ is $0$. We show that every compact composition operator $C_\phi$ on ${\cal D}_\ast$ is compact on the 
Bergman-Orlicz space ${\mathfrak B}^{\Psi_2}$ and on the Hardy-Orlicz space $H^{\Psi_2}$; in particular, $C_\phi$ is in every Schatten class $S_p$, $p > 0$, 
both on the Hardy space $H^2$ and on the Bergman space ${\mathfrak B}^2$. On the other hand, there exists a Schur function $\phi$ such that $C_\phi$ is 
compact on $H^{\Psi_2}$, but which is not even bounded on ${\cal D}_\ast$. We prove that for every $p > 0$, there exists a symbol $\phi$ such that 
$C_\phi \in S_p ({\cal D}_\ast)$, but $C_\phi \notin S_q ({\cal D}_\ast)$ for any $q < p$, that there exists another symbol $\phi$ such that 
$C_\phi \in S_q ({\cal D}_\ast)$ for every $q < p$, but $C_\phi \notin S_p ({\cal D}_\ast)$. Also, there exists a Schur function $\phi$ such that $C_\phi$ is 
compact on ${\cal D}_\ast$, but in no Schatten class $S_p ({\cal D}_\ast)$.}
\medskip

\noindent{\bf MSC 2010.} Primary: 47B33 -- Secondary: 28A12; 30C85; 31A15; 46E20; 46E22; 47B10.

\medskip

\noindent{\bf Key-words.} Bergman space -- Bergman-Orlicz space -- composition operator -- Dirichlet space -- Hardy space -- Hardy-Orlicz space -- 
logarithmic capacity -- Schatten classes 


\section{Introduction, notation and background} 

\subsection{Introduction} 

Recall that a Schur function is an analytic self-map of the open unit disk $\D$. Every Schur function $\phi$ generates a bounded composition operator 
$C_\phi$ on the Hardy space $H^2$, given by $C_\phi (f) = f \circ \phi$. Let us also introduce the set $E_\phi$ of contact points of the symbol with the unit 
circle (equipped with its normalized Haar measure $m$), namely: 
\begin{equation} 
E_\phi = \{ \e^{it} \, ; \ |\phi^\ast (\e^{it}) | = 1 \}. 
\end{equation}

In terms of $E_\phi$, a well-known necessary condition for compactness of $C_\phi$ on $H^2$ is that $m (E_\phi) = 0$. This set $E_\phi$ is otherwise more or 
less arbitrary. Indeed, it was proved in \cite{Gallardo 3} that there exist compact composition operators $C_\phi$ on $H^2$  such that the Hausdorff dimension 
of $E_\phi$ is $1$. This was generalized in \cite{ELFKELSHAYOU}: for every Lebesgue-negligible compact set $K$ of the unit circle $\T$, there is a 
Hilbert-Schmidt composition operator $C_\phi$ on $H^2$  such that  $E_\phi = K$, and in \cite{LIQUEFRODR}: 
\begin{theorem} [\cite{LIQUEFRODR}] 
For every Lebesgue-negligible compact set $K$ of the unit-circle $\T$ and every vanishing sequence $(\eps_n)$ of positive numbers, there is a  composition 
operator $C_\phi$ on $H^2$  such that $E_\phi = K$ and such that its approximation numbers satisfy $a_n (C_\phi) \leq C \, \e^{- n \, \eps_n}$.
\end{theorem} 

\medskip

We are interested here in a different Hilbert space of analytic functions, on which not every Schur function defines a bounded composition operator, namely the 
Dirichlet space $\mathcal{D}$.  Recall its definition: the Dirichlet space $\mathcal{D}$ is the space of analytic functions $f \colon \D \to \C$ such that: 
\begin{equation} 
\| f \|_{\cal D}^2 := | f (0)|^2 + \int_\D |f ' (z) |^2 \, dA (z) < + \infty \, .
\end{equation} 
If $f (z) = \sum_{n = 0}^\infty c_n z^n$, one has:
\begin{equation} 
\| f \|_{\cal D}^2 = |c_0|^2 + \sum_{n = 1}^\infty n \, |c_n|^2 \, .
\end{equation} 
Then $\| \ \|_{\cal D}$ is a norm on ${\cal D}$, making ${\cal D}$ a Hilbert space. Whereas every Schur function $\phi$ generates a bounded composition 
operator $C_\phi$ on the Hardy space $H^2$, it is no longer the case for the Dirichlet space (see \cite{McCluer-Shapiro}, Proposition 3.12, for instance). 

\medskip
In \cite{Gallardo}, the study of compact composition operators on the Dirichlet space $\mathcal{D}$ associated with a Schur function $\phi$ in connection 
with the set $E_\phi$ was initiated. In particular, it is proved there that if the composition operator $C_\phi$ is Hilbert-Schmidt on $\mathcal{D}$, then the 
logarithmic capacity ${\rm Cap}\, E_\phi$  of $E_\phi$ is $0$, but, on the other hand, there are compact composition operators on $\mathcal{D}$ for which 
this capacity is positive. The optimality of this theorem was later proved in \cite{ELFKELSHAYOU} under the following form: 
\begin{theorem} [O. El-Fallah, K. Kellay, M. Shabankhah, H. Youssfi] \hfill \\ 
For every compact set $K$ of the unit circle $\T$ with  logarithmic capacity ${\rm Cap}\, K$ equal to $0$, there exits a Hilbert-Schmidt composition operator 
$C_\phi$ on $\mathcal{D}$  such that  $E_\phi = K$. 
\end{theorem} 
\par

In this paper, we shall improve on this last result. We prove in Section~\ref{section capa} (Theorem~\ref{main th}) that for every compact set 
$K \subseteq \partial{\D}$ of logarithmic capacity ${\rm Cap} \, K = 0$, there exists a Schur function $\phi \in A (\D) \cap {\cal D}_\ast$ such that 
the composition operator $C_\phi$ is in all Schatten classes $S_p ({\cal D}_\ast)$, $p > 0$, and for which $E_\phi = K$ (and moreover 
$E_\phi = \{\e^{it} \, ; \ \phi (\e^{it}) = 1\}$). On the other hand, in Section~\ref{section boundedness}, we show (Theorem~\ref{capa nulle}) that for 
every bounded composition operator $C_\phi$ on  ${\cal D}_\ast$ and every $\xi \in \partial \D$, the logarithmic capacity of 
$E_\phi (\xi) = \{ \e^{it} \, ; \ \phi (\e^{it}) = \xi \}$ is $0$. \par

In link with Hardy and Bergman spaces, we prove, in Section~\ref{section boundedness} yet, that every compact composition operator $C_\phi$ on 
${\cal D}_\ast$ is compact on the Bergman-Orlicz space ${\mathfrak B}^{\Psi_2}$ and on the Hardy-Orlicz space $H^{\Psi_2}$. In particular, $C_\phi$ is in 
every Schatten class $S_p$, $p > 0$, both on the Hardy space $H^2$ and on the Bergman space ${\mathfrak B}^2$ (Theorem~\ref{compact Dirichlet and 
Schatten Hardy}).  However, there exists a Schur function $\phi$ such that $C_\phi$ is compact on $H^{\Psi_2}$, but which is not even bounded on 
${\cal D}_\ast$ (Theorem~\ref{not bounded}). \par

In Section~\ref{section Schatten}, we give a characterization of the membership of composition operators in the Schatten classes $S_p ({\cal D}_\ast)$, $p > 0$ 
(actually in $S_p ({\cal D}_{\alpha, \ast} ) $, where ${\cal D}_{\alpha, \ast}$ is the weighted Dirichlet space). We deduce that for every $p > 0$, there exists 
a symbol $\phi$ such that $C_\phi \in S_p ({\cal D}_\ast)$, but $C_\phi \notin S_q ({\cal D}_\ast)$ for any $q < p$, and that there exists another symbol 
$\phi$ such that $C_\phi \in S_q ({\cal D}_\ast)$ for every $q < p$, but $C_\phi \notin S_p ({\cal D}_\ast)$ (Theorem~\ref{separation Schatten}). We also 
show that there exists a Schur function $\phi$ such that $C_\phi$ is compact on ${\cal D}_\ast$, but in no Schatten class $S_p ({\cal D}_\ast)$ 
(Theorem~\ref{no Schatten}). 

\goodbreak 
 
\subsection{Notation and background.} 

We denote by $\D$ the unit open disk of the complex plane and by $\T = \partial \D$ the unit circle. $A$ is the normalized area measure $dx \, dy / \pi$ of 
$\D$ and $m$ the normalized Lebesgue measure $dt/ 2 \pi$ on $\T$. \par

As said before, a Schur function is an analytic self-map of $\D$ and the associated composition operator is defined, formally, by $C_\phi (f) = f \circ \phi$. 
The function $\phi$ is called the symbol of $C_\phi$.
\par

The Dirichlet space ${\cal D}$ is defined above. We shall actually work, for convenience, with its subspace ${\cal D}_\ast$ of functions $f \in {\cal D}$ such 
that $f (0) = 0$. In this paper, we call ${\cal D}_\ast$ the \emph{Dirichlet space}. \par 
\smallskip 

An orthonormal basis of $\mathcal{D}_\ast$ is formed by $e_n (z) = z^n / \sqrt{n}$, $n \geq 1$. The reproducing kernel on $\mathcal{D}_\ast$, defined by 
$f (a) =\langle f , K_a \rangle$ for every $f \in {\cal D}_\ast$, is given by $K_a (z) = \sum_{n = 1}^\infty \overline{e_n (a)} \, e_n (z)$, so that:
\begin{equation} \label{repro} 
K_a (z) = \log \frac{1}{1 - \overline{a} z} \, \cdot 
\end{equation}

Compactness of composition operators on ${\cal D}$ was characterized in terms of Carleson measure by D. Stegenga (\cite{Stegenga}) and by B. McCluer and 
J. Shapiro in terms of angular derivative (\cite{McCluer-Shapiro}). Another characterization, more useful for us here, was given by 
N. Zorboska (\cite{Zor}, page~2020): for $\phi \in {\cal D}$, $C_\phi$ is bounded on ${\cal D}$ if and only:
\begin{equation} \label{bounded} 
\sup_{h \in (0, 2)} \sup_{|\xi| = 1} \frac{1}{A [ W (\xi, h)]} \int_{W (\xi, h)} n_\phi (w) \, dA (w) < \infty \, , 
\end{equation} 
where $W (\xi, h) = \{ w \in \D \, ; \ 1 - |w| \leq h \text{ and } |\arg (w \bar{\xi}) | \leq \pi h\}$ is the Carleson window of size $h \in (0, 2)$ center at 
$\xi \in \T$ and $n_\phi$ is the counting function of $\phi$:
\begin{equation} 
\qquad n_\phi (w) = \sum_{\phi (z) = w} 1 \ , \qquad w \in \phi (\D) \, , 
\end{equation} 
(we set $n_\phi (w) = 0$ for $w \in \D \setminus \phi(\D)$). In particular, every Schur function with bounded valence defines a bounded composition 
operator on ${\cal D}$. \par

Moreover, $C_\phi$ is compact if and only if:
\begin{equation} \label{compact} 
\sup_{|\xi| = 1} \frac{1}{A [ W (\xi, h)]} \int_{W (\xi, h)} n_\phi (w) \, dA (w) \mathop{\longrightarrow}_{h \to 0} 0 \, .
\end{equation} 
\smallskip 

For further informations on the Dirichlet space, one may consult the two surveys \cite{survey} and \cite{Ross}, for example. 
\goodbreak 


\subsubsection{Logarithmic capacity} \label{rappels capa}

The notion of logarithmic capacity is tied to the study of the Dirichlet space by the following seminal and sharp result of Beurling (\cite{BEU}; see also 
\cite{KAHSAL}).
\begin{theorem} [Beurling] 
For every function $f (z) = \sum_{n = 0}^\infty c_n z^n \in \mathcal{D}$, there exists a set $E \subseteq \partial{\D}$, with logarithmic 
capacity $0$, such that, if $t\in \T \setminus E$, then the radial limit  $f^\ast (\e^{it}) := \lim_{r \to 1^{-}} f (r \e^{it})$ exists (in $\C$). Moreover, the result 
is optimal: if a compact set $E \subseteq \T$ has zero logarithmic capacity, there exists $f (z) = \sum_{n = 0}^\infty c_n z^n \in \mathcal{D}$ such that 
$f^\ast (\e^{it})$ does not exist on $E$. 
\end{theorem} 
\medskip

Let us recall some definitions (see \cite{KAHSAL}, Chapitre~III, \cite{Conway}, Chapter~21, \S~7, or \cite{Ross}, Section~4, for example). \par
\smallskip

Let $\mu$ be a probability measure supported by a compact subset $K$ of $\T$. The \emph{potential} $U_{\mu}$ of $\mu$ is defined, for every $z \in \C$, by:
\begin{displaymath}
U_{\mu} (z) = \int_K \log \frac{\e}{|z - w|} \, d\mu (w) \, .
\end{displaymath} 
The \emph{energy} $I_\mu$ of $\mu$ is defined by:
\begin{displaymath}
I_\mu = \int_K U_{\mu} (z) \, d\mu (z) = \iint_{K \times K} \log \frac{\e}{|z - w|} \, d\mu (w) \, d\mu (z) \,.
\end{displaymath} 
The \emph{logarithmic capacity} of a Borel set $E \subseteq \T$ is:
\begin{displaymath} 
{\rm Cap} \, E = \sup\nolimits_\mu \e^{- I_\mu} \,, 
\end{displaymath} 
where the supremum is over all Borel probability measures $\mu$ with compact support contained in $E$ . Hence $E$ is of logarithmic capacity $0$ (which is the 
case we are interested in) if  and only if $I_\mu = \infty$ for all probability measures compactly carried by $E$. The fact that ${\rm Cap} \, E = 0$ implies 
that $E$ has null Lebesgue measure (\cite{KAHSAL}, Chapitre~III, Th\'eor\`eme~I) (hence ${\rm Cap} \, E > 0$ if $E$ is a non-void open subset of $\T$), but 
the converse is wrong, as shown by Cantor's middle-third set ${\mathfrak C}$. A compact set $K$ such that ${\rm Cap}\, K = 0$ is totally disconnected 
(\cite{Conway}, Corollary~21.7.7). \par
\smallskip

If $E$ is a compact set with ${\rm Cap} \, E > 0$, there is a unique probability measure compactly carried by $E$ that minimizes the energy 
$I_\mu$ (\cite{Conway}, Theorem~21.10.2, or \cite{KAHSAL}, Chapitre~III, Proposition~4). Such a measure is called the \emph{equilibrium measure} 
of $E$. 
\par 
If $\mu$ is the equilibrium measure of the compact set $K$, we have Frostman's Theorem (\cite{Conway}, Theorem~21.7.12, or \cite{KAHSAL}, Chapitre~III, 
Proposition~5 and Proposition~6): $U_\mu (z) \leq I_\mu$ for every $z \in \C$ and 
\begin{equation} \label{almosteverywhere} 
\qquad U_{\mu} (z) = I_\mu \quad \text{for almost all } z \in K \, .  
\end{equation} 

\medskip

Suppose that the compact set $K$ has zero logarithmic capacity. For $\eps > 0$, let  $K_\eps = \{ z \in \T \, ; \ {\rm dist} \, (z, K) \leq \eps \}$,  
$\mu_\eps$ its equilibrium measure, and $I_{\mu_\eps}$ its energy. Then (\cite{Conway}, Proposition~21.7.15):
\begin{equation} \label{explose}  
\lim_{\eps \to 0} I_{\mu_\eps} = \infty \, .
\end{equation} 

\goodbreak 


\section{Bounded and compact composition operators}  \label{section boundedness}

In \cite{Gallardo},  E. A. Gallardo-Guti\'errez and M. J. Gonz\'alez showed that for every Hilbert-Schmidt composition operator $C_\phi$ on ${\cal D}_\ast$, 
the logarithmic capacity of the set $E_\phi = \{ \e^{i\theta} \in \partial \D \, ; \ |\phi (\e^{i \theta}) | = 1 \}$ is zero. On the other hand, they showed that 
there are compact composition operators on ${\cal D}_\ast$ for which $E_\phi$ has positive logarithmic capacity. We shall see that if we replace $|\phi|$ by 
$\phi$ in the definition of $E_\phi$, the result is very different.

\begin{theorem} \label{capa nulle} 
For every bounded composition operator $C_\phi$ on  ${\cal D}_\ast$ and every $\xi \in \partial \D$, the logarithmic capacity of 
$E_\phi (\xi) = \{ \e^{it} \, ; \ \phi (\e^{it}) = \xi \}$ is $0$. 
\end{theorem} 

We first state the following characterization of Hilbert-Schmidt composition operators on ${\cal D}_\ast$. This result is stated in \cite{Gallardo}, but not 
entirely proved.

\begin{lemma}
Let $\phi \in {\cal D}_\ast$ be an analytic self-map of $\D$. Then $C_\phi$ is Hilbert-Schmidt on ${\cal D}_\ast$ if and only if 
\begin{equation} \label{int HS finie} 
\int_\D \frac{|\phi ' (z)|^2}{(1 - |\phi (z) |^2)^2} \, dA (z) < \infty \, .
\end{equation} 
\end{lemma}

\noindent{\bf Proof.} Let $e_n (z) = z^n / \sqrt n$; then $(e_n)_{n \geq 1}$ is an orthonormal basis of ${\cal D}_\ast$ and 
\begin{displaymath} 
\sum_{n = 1}^\infty \| C_\phi (e_n) \|^2 = \sum_{n = 1}^\infty \frac{\| \phi^n\|^2}{n} 
= \int_\D \frac{|\phi ' (z)|^2}{(1 - |\phi (z) |^2)^2} \, dA (z) \,. 
\end{displaymath} 
Hence \eqref{int HS finie} is satisfied if $C_\phi$ is Hilbert-Schmidt. To get the converse, we need to show that \eqref{int HS finie} implies that $C_\phi$ is 
bounded on ${\cal D}_\ast$. Let $f \in {\cal D}_\ast$ and write $f (z) = \sum_{n = 1}^\infty c_n z^n$. Then 
$C_\phi f = \sum_{n = 1}^\infty c_n \phi^n$ and
\begin{align*}
\| C_\phi f \| 
& \leq \sum_{n = 1}^\infty |c_n|\, \| \phi^n \|  
\leq \bigg(\sum_{n = 1}^\infty n \, |c_n|^2 \bigg)^{1/2}  \bigg( \sum_{n = 1}^\infty \frac{\| \phi^n\|^2}{n} \bigg)^{1/2} \\
& = \bigg(\int_\D \frac{|\phi ' (z)|^2}{(1 - |\phi (z) |^2)^2} \, dA (z) \bigg)^{1/2} \| f \| \, .
\end{align*}
Then \eqref{int HS finie} implies that $C_\phi$ is Hilbert-Schmidt. 
\qed

\bigskip

Now Theorem~\ref{capa nulle} will follow from the next proposition.
\begin{proposition} \label{existence HS} 
There exists an analytic self-map $\sigma$ of $\D$, belonging to ${\cal D}_\ast$ and to the disk algebra $A (\D)$, such that $\sigma (1) = 1$ and 
$|\sigma (\xi) | < 1$ for $\xi \in \partial \D \setminus \{1\}$ and such that the associated composition operator $C_\sigma$ is Hilbert-Schmidt on 
${\cal D}_\ast$. 
\end{proposition}

Taking this proposition for granted for a while, we can prove the theorem.\par
\medskip

\noindent{\bf Proof of Theorem~\ref{capa nulle}.} Making a rotation, we may, and do, assume that $\xi = 1$. Then, if $\sigma$ is the map of 
Proposition~\ref{existence HS},  $C_\phi C_\sigma = C_{\sigma \circ \phi}$ is Hilbert-Schmidt. By \cite{Gallardo}, the set 
$E_{\sigma \circ \phi}$ has zero logarithmic capacity. But $\sigma$ has modulus $1$ only at $1$; hence $\e^{i \theta} \in E_{\sigma \circ \phi}$ if and only 
if $\e^{i \theta} \in E_\phi (1)$. 
\qed

\bigskip

To prove Proposition~\ref{existence HS}, it will be convenient to use the following criteria, where $\phi_a (z) = \frac{z - a}{1 -\bar{a} z} \, \cdot$ 

\begin{lemma} \label{CS HS} 
Let $f \in {\cal D}$ such that $\Re f \geq 1$. Then if $\sigma = \phi_a \circ \e^{- 1 / f } $, where $a = \e^{- 1 / f (0)}$, the composition operator $C_\sigma$ is 
Hilbert-Schmidt on ${\cal D}_\ast$. 
\end{lemma}

\noindent{\bf Proof.} Let $\sigma_0 = \e^{- 1 / f }$. If $u = \Re f $ and $v = \Im f $, one has:
\begin{displaymath} 
|\sigma_0|^2 = \exp \Big( - \frac{2 u}{u^2 + v^2} \Big) \quad \text{and} \quad  
|\sigma_0 ' |^2 = \frac{{u '}^2 + {v '}^2}{(u^2 + v^2)^2} \, \exp \Big( - \frac{2 u}{u^2 + v^2} \Big) \, \cdot
\end{displaymath} 
Then $|\sigma_0| < 1$ and so $\sigma_0$ is a self-map of $\D$. Since $u \geq 1 > 0$, one has 
$|\sigma_0 '|^2 \leq ({u '}^2 + {v '}^2) / (u^2 + v^2)^2 \leq {u '}^2 + {v '}^2 = |f ' |^2$; hence $\sigma_0 \in {\cal D}$. \par
For $0 \leq x \leq 2$, one has $1 - \e^{- x} \geq x / 4$. Therefore, since $u \geq 1$ implies $2u / (u^2 + v^2) \leq 2 / u \leq 2$, one has:
\begin{displaymath} 
1 - |\sigma_0|^2 \geq \frac{u}{2 (u^2 + v^2)} \, \cdot
\end{displaymath} 
It follows that:
\begin{displaymath} 
\frac{|\sigma_0 '|^2}{(1 - |\sigma_0|^2)^2} \leq \frac{{u '}^2 + {v '}^2}{(u^2 + v^2)^2} \, \frac{4  (u^2 + v^2)^2}{u^2} 
\leq 4 ({u '}^2 + {v '}^2) = 4 |f '|^2 \,.
\end{displaymath} 
Since $f \in {\cal D}$, $|f ' |^2$ has a finite integral and therefore \eqref{int HS finie} is satisfied. It follows that $C_{\sigma_0}$ is Hilbert-Schmidt on 
${\cal D}$ and hence $C_\sigma = C_{\sigma_0} \circ C_{\phi_a}$ is Hilbert-Schmidt on ${\cal D}_\ast$, since $\sigma (0) = 0$. 
\qed

\bigskip 

\noindent{\bf Proof of Proposition~\ref{existence HS}.} Let $\Omega$ be the domain defined by:
\begin{displaymath} 
\Omega = \{ z \in \C \, ; \ \Re z > 1 \text{ and } |\Im z| < 1/ (\Re z)^2 \} \, .
\end{displaymath} 
Let $f$ be a conformal map from $\D$ onto $\Omega$ such that $f (1) = \infty$. Since $A (\Omega) < \infty$, we have $f \in {\cal D}$. By 
Lemma~\ref{CS HS}, the function $\sigma = \e^{- 1 / f}$ has the required properties. 
\qed 

\bigskip 

For the next result, recall that an Orlicz function $\Psi$ is a nondecreasing convex function such that $\Psi (0) = 0$ and $\Psi (x) / x \to \infty$ as $x$ goes to 
infinity. We refer to \cite{memoirs} for the definition of Hardy-Orlicz and Bergman-Orlicz spaces. In the following result, one set $\Psi_2 (x) = \exp (x^2) - 1$. 

\goodbreak
\begin{theorem} \label{compact Dirichlet and Schatten Hardy} 
Every compact composition operator $C_\phi$ on ${\cal D}_\ast$ is compact on the Bergman-Orlicz space ${\mathfrak B}^{\Psi_2}$ and on the Hardy-Orlicz 
space $H^{\Psi_2}$. In particular, $C_\phi$ is in every Schatten class $S_p$, $p > 0$, both on the Hardy space $H^2$ and on the Bergman space 
${\mathfrak B}^2$. 
\end{theorem} 

\noindent{\bf Proof.} Consider the normalized reproducing kernels $\tilde K_a = K_a /  \| K_a\|$, $a \in \D$. When $|a|$ goes to $1$, they tends to $0$ 
uniformly on compact sets of $\D$; hence $\| C_\phi^\ast ( \tilde K_a) \|$ tends to $0$, by compactness of the adjoint operator $C_\phi^\ast$. But 
$C_\phi^\ast (K_a) = K_{\phi (a)}$ and $\| K_a \|^2 = \langle K_a , K_a \rangle = \log \frac{1}{1 - |a|^2}$\raise 1pt \hbox{,} so we get:
\begin{equation} \label{CN compacite} 
\lim_{|a| \to 1} \frac{ \log \frac{1}{1 - |\phi (a)|^2} } {\log \frac{1}{1 - |a|^2} } = 0 \,. 
\end{equation}
This condition means that $C_\phi$ is compact on the Bergman-Orlicz space ${\mathfrak B}^{\Psi_2}$ (\cite{memoirs}, page~69) and implies that 
$C_\phi$ is in all Schatten classes $S_p ({\mathfrak B}^2)$, $p > 0$ (\cite{LLQR Bergman}). \par 
\medskip

In the same way, it suffices to show that $C_\phi$ is compact on $H^{\Psi_2}$, because that implies that $C_\phi$ is in all Schatten classes $S_p (H^2)$ 
(\cite{LLQR 2009}, Theorem~5.2). \par
Compactness of $C_\phi$ on $H^{\Psi}$ is equivalent to say (\cite{memoirs}, Theorem~4.18) that:
\begin{align*} 
\rho_\phi (h) := \sup_{|\xi| = 1} m \big( \{\e^{it} \, ; \ \phi (\e^{it}) & \in W (\xi, h) \} \big) \\ 
& = o_{h \to 0} \, \bigg[ \frac{1}{\Psi \big( A \Psi^{- 1} (1 / h) \big) } \bigg] 
\quad \text{for every } A > 0 \,. 
\end{align*} 
When $\Psi = \Psi_2$, this means that $\rho_\phi (h) = o\, (h^A)$ for every $A > 0$. Now, by \cite{Nevanlinna}, Theorem~4.2, this is also equivalent to say 
that:
\begin{equation} \label{petito} 
\sup_{|\xi| = 1} \frac{1}{A [W (\xi, h)]} \int_{W (\xi, h)} N_\phi (w) \, dA (w) = o\, (h^A) \quad \text{for every } A > 0 \, , 
\end{equation} 
where $N_\phi$ is the Nevanlinna counting function of $\phi$:
\begin{equation} \label{Nevanlinna} 
\qquad N_\phi (w) = \sum_{\phi (z) = w} (1 - |z|^2) \, , \qquad w \in \phi (\D) \,, 
\end{equation} 
and $N_\phi (w) = 0$ otherwise. \par 

But \eqref{CN compacite} is equivalent to the fact that for every $\eps > 0$ there exists $\delta_\eps > 0$ such that:
\begin{equation} \label{necess}  
\qquad 1 -| \phi (z) | \geq \delta_\eps (1 - | z |)^\eps \, , \qquad  \forall z \in \D \,.  
\end{equation}
Since $\phi (0) = 0$, we have $|\phi (z) | \leq |z|$, by Schwarz's lemma; hence one has 
$N_\phi (w) \leq 2 \delta_\eps^{- 1} (1 - |w|)^{1/ \eps} n_\phi (w)$. It follows that (since $1 - |w| \leq h$ for $w \in W (\xi, h)$):
\begin{displaymath} 
\frac{1}{A [W (\xi, h)]} \int_{W (\xi, h)} N_\phi (w) \, dA (w) 
\leq 2 \delta_\eps^{ - 1} h^{1/ \eps} \frac{1}{A [W (\xi, h)]} \int_{W (\xi, h)} n_\phi (w) \, dA (w) \, , 
\end{displaymath} 
which is $o\, (h^{1/\eps})$, uniformly for $|\xi| = 1$, by \eqref{compact}. 
\qed

\bigskip

\noindent{\bf Remarks.} 1) One may argue that compactness of $C_\phi$ on $H^{\Psi_2}$ implies its compactness on ${\mathfrak B}^{\Psi_2}$ 
(\cite{LLQR Bergman}, Proposition~4.1, or \cite{Racsam}, Theorem~9). One may also use the forthcoming 
Corollary~\ref{SchattenDirichletalphaimpliqueSchattenDirichletbeta} saying that $C_\phi \in S_p (H^2)$ implies that $C_\phi \in S_p ({\mathfrak B}^2)$. 
\par 
2) To show the compactness of $C_\phi$ on $H^{\Psi_2}$, we used its compactness on ${\cal D}_\ast$ twice. However, due to the fact that $\eps > 0$ is 
arbitrary, we may replace $o\, (h^{1/ \eps})$ by $O\, (h^{1/ \eps})$; hence to end the proof, we only have to use \eqref{bounded}, {\it i.e.} the boundedness of 
$C_\phi$ on ${\cal D}_\ast$, instead of \eqref{compact}. \par

Note that \eqref{CN compacite} does not suffice to have compactness on $H^{\Psi_2}$ (in \cite{memoirs}, Proposition~5.5, we construct a Blaschke product 
satisfying \eqref{CN compacite}). 

\bigskip

In the opposite direction, we have the following result. 

\goodbreak 

\begin{theorem} \label{not bounded} 
There exists a Schur function $\phi$ such that $C_\phi$ is compact on $H^{\Psi_2}$, but which is not even bounded on ${\cal D}_\ast$. 
\end{theorem}

To prove this theorem, we first begin with the following key lemma.
\goodbreak
\begin{lemma} \label{lemme sans barriere}
There exists a constant $\kappa_1 > 0$ such that for any $f \in {\cal H} (\D)$ having radial limits $f^\ast$ a.e. and which satisfies, for some $\alpha \in \R$:
\begin{equation} \label{Luis 2}
\left\{
\begin{array} {l}
\Im f (0) < \alpha \qquad \text{and} \\ 
f (\D) \subseteq \{z \in \C \, ; \  0 < \Re z < \pi \} \cup \{z \in \C \, ; \  \Im z < \alpha \} ,
\end{array}
\right. 
\end{equation} 
we have, for all $y \geq \alpha$:
\begin{displaymath} 
m \big( \{ z \in \T \, ; \  \Im [ f^\ast (z) ] \geq y \} \big) \leq \kappa_1 \e^{\alpha - y} .
\end{displaymath} 
\end{lemma}

\noindent{\bf Proof.} Suppose that $f$ satisfies \eqref{Luis 2}, and define $f_1 (z) = - i f (z) + \frac{\pi}{2} i - \alpha$. Then either $\Re [ f_1 (z)] < 0$, 
or $- \frac{\pi}{2} < \Im [f_1 (z)] < \frac{\pi}{2}$ for every $z \in \D$. Therefore, defining $h (z) = 1 + \exp [f_1 (z)]$, we have $h \colon \D \to \H$, 
that is $\Re [h (z)] > 0$ for every $z \in \D$. \par
Finally define $h_1 (z) = h (z) - i \, \Im [h (0)]$. Then $h_1 \colon \D \to \H$ and $h_1 (0) \in \R$ (and so $h_1 (0) > 0$). Kolmogorov's inequality yields that, 
for some absolute constant $C_1$, one has, for every $\lambda > 0$:
\begin{equation} \label{Kolmogorov}
m \big( \{ z \in \T \, ; \ | h_1^\ast (z) | \geq \lambda \} \big) \leq C_1 \, \frac{h_1 (0)}{\lambda} \, \cdot
\end{equation} 
Observe that, since $\Im [ f (0)] < \alpha$, we have $\Re [ f_1 (0)] < 0$, and then:
\begin{equation} \label{Luis 3}
| \Im [h (0)] | < 1 \qquad \text{and} \qquad h_1 (0) = \Re [h (0)]  < 2.
\end{equation} 
Suppose now that, for $y > \alpha$ and $z \in \D$, we have $\Im [f (z)] > y$; then $\exp [f_1 (z)] \in \H$, and 
$| h (z) | \geq | \exp [ f_1 (z)] | > \e^{y - \alpha}$. Taking radial limits we get, up to a set of null Lebesgue-measure:
\begin{displaymath} 
\{ z \in \T \, ; \ \Im [ f^\ast (z)] \geq y \} \subseteq \{ z \in \T \, ; \ | h^\ast (z) | \geq \e^{y - \alpha} \} .
\end{displaymath} 

We consider two cases: $\e^{y - \alpha} \geq 2$ and $\e^{y - \alpha} < 2$. When $\e^{y - \alpha} \geq 2$, then $| h^\ast ( z) | \geq \e^{y - \alpha}$ 
yields:
\begin{displaymath} 
| h_1^\ast (z) | \geq \e^{y - \alpha} - |\Im [h (0)] | > \e^{y - \alpha} - 1 \geq \frac{1}{2}\, \e^{y - \alpha} , 
\end{displaymath} 
by the first part of \eqref{Luis 3}. Then, using \eqref{Kolmogorov} and the second part of \eqref{Luis 3}, we have:
\begin{align*} 
m \big( \{ z \in \T \, ; \ \Im [ f^\ast (z)] \geq y \} \big) 
& \leq m \big( \{ z \in \T \, ; \ | h_1^\ast (z) | > (1/2) \, \e^{y - \alpha} \} \big) \\
& \leq \frac{2 \, C_1 h_1 (0)}{\e^{y - \alpha}} \leq \frac{4 C_1}{\e^{y - \alpha}} \, \raise 1pt \hbox{,} 
\end{align*} 
and, in this case, the lemma is proved, if one takes $\kappa_1 \geq 4 C_1$. \par

When $\e^{y - \alpha} < 2$, then $\e^{\alpha - y} > 1/2$, and, because:
\begin{displaymath} 
m \big( \{ z \in \T \, ; \ \Im [ f^\ast (z)] \geq y \} \big) \leq 1 < \kappa_1 \e^{\alpha - y} ,
\end{displaymath} 
since $\kappa_1 > 2$, the lemma is proved. 
\qed 
\medskip

Now, we give a general construction of Schur functions with suitable properties. 

\begin{proposition} \label{general construction} 
Let ${\mathfrak g} \colon (0, \infty) \to (0, \infty)$ be a continuous non-increasing function such that:
\begin{displaymath} 
\lim_{t \to 0^+} {\mathfrak g} (t) = + \infty , \qquad \text{and} \qquad \lim_{t \to + \infty} {\mathfrak g} (t) = 0 . 
\end{displaymath} 
Let ${\mathfrak h} \colon (0, \infty) \to (0, \infty]$ be a lower semicontinuous function such that 
$M := \sup\{{\mathfrak h}(t) \, ; \ t \ge \pi \} < + \infty \, $ and consider the simply connected domain:
\begin{displaymath} 
\Omega = \{ x + i y \,  ; \ x \in (0, \infty) \quad \text{and} \quad  {\mathfrak g} (x) < y < {\mathfrak g} (x) + {\mathfrak h}(x) \} \, . 
\end{displaymath}
Let ${\mathfrak f}  \colon \overline{\D} \to \overline{\Omega} \cup \{\infty\}$ be a conformal mapping from $\D$ onto $\Omega$  
such that ${\mathfrak f} (0) = \pi +  i ({\mathfrak g}(\pi) +{\mathfrak h}(\pi) /2 )$.
\par\smallskip 

Then the symbol $\phi \colon \D \to \D$ defined by $\phi (z) = \exp [ - {\mathfrak f} (z)]$, for every $z \in \D$, satisfies, for some $\eps_0$,$ k_0 > 0$: \par
\smallskip 

1) For all $h \in (0, \eps_0)$:
\begin{equation} \label{Luis 1} 
m ( \{z \in \T \, ; \  |\phi^\ast (z)| > 1 - h \}) \leq k_0 \exp \big(- {\mathfrak g} (2h) \big) \, .
\end{equation} 
\par

2) Assume that, for some $r \in (0,\infty]$ and integers $0 \leq n < N \leq \infty$, one has 
$\{{\mathfrak h}(t) \, ; \  t \le r \}\subseteq  \,  ( 2 \, n \pi, 2 N \pi]$. Then, for all $z\in\D$, such that $|z|> \e^{- r}$, we have  $n \le n_\phi (z) \le N$. 
\par
In particular, $\{z \in \D \, ; \ |z|> \e^{- r}\} \subseteq \phi (\D) \subseteq \D \setminus \{0\}$, when $n \ge 1$.
\end{proposition}

\noindent{\bf Remarks.} \par\smallskip 

1. When $N = 1$, the map $\phi$ is univalent. \par\smallskip 

2. When $r = \infty$ and $n \ge 1$, we have $\phi (\D) = \D\setminus\{0\}$. \par\smallskip 

3. With ${\mathfrak g} (t) =1/ t$, the operator $C_\phi$ is compact on $H^{\Psi_2}$, therefore belongs to all Schatten classes $S_p (H^2)$, $p > 0$. 
\par\smallskip 

4. When $N < \infty$, the operator $C_\phi$ is bounded on the Dirichlet space. \par\smallskip 

5. When $n \geq 1$, the operator $C_\phi$ is not compact on the Dirichlet space (since the averages on the windows  of the function $n_\phi$ cannot uniformly 
vanish).
\par\medskip

\noindent{\bf Proof of Proposition~\ref{general construction}.} We shall apply Lemma~\ref{lemme sans barriere} with $\alpha = M + {\mathfrak g} (\pi)$. 
\par\smallskip 

Suppose that, for $z \in \T$ and $0 < h < 1$, we have $|\phi^\ast ( z) | > 1 - h$. Then, if $h$ is small enough, 
\begin{displaymath} 
\e^{- 2 h} < 1 - h < |\phi^\ast (z)| = \exp \big( - \Re [f^\ast (z)] \big) ,
\end{displaymath} 
and therefore $2h > \Re [ {\mathfrak f}^\ast (z)]$. But observe that ${\mathfrak f}^\ast (z) \in \overline{\Omega} \cup \{\infty\}$, and so, if 
$2 h > \Re [f^\ast (z)]$, we necessarily have $\Im [{\mathfrak f}^\ast (z)] \geq {\mathfrak g} (2h)$. Again, if $h$ is small enough, we have 
$y = {\mathfrak g} (2h) > \alpha$, and may apply the lemma to obtain:
\begin{displaymath} 
m \big( \{ z \in \T \, ; \ |\phi^\ast (z) | > 1 - h \} \big) \leq m \big( \{ z \in \T \, ; \ \Im [ {\mathfrak f}^\ast (z)] \geq {\mathfrak g} (2h) \} \big) 
\leq \kappa_1 \e^{\alpha - {\mathfrak g} (2h)} .
\end{displaymath} 
We get \eqref{Luis 1}. \par

On the other hand, let $Z\in\D$ such that $|Z|> \e^{- r}$, we can write $Z = \e^{- x} \e^{i \theta}$ with $x < r$. We can find $\theta_i's$ such that  
${\mathfrak g} (x)< \theta_1 < \ldots < \theta_s < {\mathfrak g} (x) + {\mathfrak h} (x)$ and  
$\theta_j \equiv \theta [2\pi]$ with $n \le s \le N$. For each $j$, there exists a unique $z_j \in \D$, such that $\Re {\mathfrak f} (z_j) = x$ and 
$\Im {\mathfrak f} (z_j) = \theta_j$; hence $\phi (z_j) = Z$. Moreover no other $z \in \D$ can satisfy $\phi (z) = Z$. Hence $n_\phi(Z) = s$.
\qed
\medskip

\noindent{\bf Proof of Theorem~\ref{not bounded}.} As said before, if one takes ${\mathfrak g} (t) = 1/ t$ in Proposition~\ref{general construction}, then 
$C_\phi$ is compact on $H^{\Psi_2}$ and hence is in all Schatten classes $S_p (H^2)$, $p > 0$. On the other hand, if one choose also 
${\mathfrak h} (t) = 1/ t$, then, for every $r > 0$, $\{ {\mathfrak h} (t) \, ; \ t \leq r \} = [1/ r , \infty)$ and for $|z | > \e^{- r}$, we get that 
$n_\phi (z) \geq [1/ (2 \pi  r) ]$ (the integer part of $1 / (2 \pi r)$). It follows that, for some constant $c > 0$, one has, with $\e^{- r} = 1 - h$:
\begin{displaymath} 
\frac{1}{A [W (\xi, h)]} \int_{W (\xi, h)} n_\phi (z) \, dA (z) \geq c\, \frac{1}{\log [1 / (1 - h)]} \mathop{\longrightarrow}_{h \to 0}  \infty \, . 
\end{displaymath} 
Therefore, $C_\phi$ is not bounded on ${\cal D}_\ast$, by \eqref{bounded}. \qed

\bigskip

\noindent{\bf Remarks.}  1. Actually, as we may take ${\mathfrak g}$ growing as we wish, the proof shows, using \cite{memoirs}, Theorem~4.18, that for every 
Orlicz function $\Psi$, one can find a Schur function $\phi$ such that $C_\phi$ is not bounded on ${\cal D}_\ast$, though compact on the Hardy-Orlicz space 
$H^\Psi$. \par 

2. This construction also allows to produce a \emph{univalent} map $\phi$, with an arbitrary small Carleson function 
$\rho_\phi (h) = \sup_{|\xi| = 1} m \big(\{ \e^{it} \, ; \ \phi^\ast (\e^{it}) \in W (\xi, h) \} \big)$, and  such that $C_\phi$ is not compact on the Dirichlet 
space (note we cannot replace ``compact'' by ``bounded'' since any Schur function with a bounded valence is bounded on the Dirichlet space). \par
Indeed, take ${\mathfrak h}(t) = 2\pi$ and ${\mathfrak g}$ be ${\cal C}^1$: ${\mathfrak g} (t) = 1/t$ for instance. We have $N = 1$ and so $\phi$ is univalent.
Now it suffices to notice that the range of the curve
\begin{displaymath} 
\Gamma = \big\{ \e^{- x -i {\mathfrak g} (x)} \, ; \  x \in (0,\infty) \big\} 
= \big\{\big( t \cos (1/\ln (t)) , t \sin (1/ \ln (t) )\big) \, ; \  t \in (0,1) \big\} \subseteq \D
\end{displaymath} 
has a null area measure. The range of $\phi$ is $\D \setminus (\Gamma \cup \{0\})$ and for each $w \notin \Gamma$, we have $n_\phi(w) = 1$
Then, for $h \in (0,1)$, we have:
\begin{align*} 
\frac{1}{h^2} \int_{W (1, h)} n_\phi (w) \, dA (w) 
& = \frac{1}{h^2} \int_{W (1, h) \setminus \Gamma} dA(w)
= \frac{1}{h^2} \, A [W (1, h) \setminus \Gamma] \\ 
& = \frac{1}{h^2} \, A [W(1, h)] \approx 1 \, ,
\end{align*} 
and so $C_\phi$ in not compact on ${\cal D}_\ast$, by \eqref{compact}. \qed  
\goodbreak 


\section {Composition operators in Schatten classes}  \label{section Schatten}

\subsection {Characterization} 

In this section, we give a characterization of the membership in the Schatten classes of composition operators on ${\cal D}_\ast$. This characterization 
will be deduced from Luecking's one for composition operators on the Bergman space. Actually, we shall give it for weighted Dirichlet spaces 
${\cal D}_{\alpha, \ast}$. Boundedness and compactness has been characterized by B. McCluer and J. Shapiro in \cite{McCluer-Shapiro} and, in other terms, 
by N. Zorboska in \cite{Zor}. \par
\medskip

Recall that for $\alpha > - 1$, the weighted Dirichlet space ${\cal D}_\alpha$ is the space of analytic functions $f \colon \D \to \C$ such that 
\begin{equation}
\int_\D |f ' (z)|^2 \, (1 - |z|^2)^\alpha \, dA (z) < \infty\, .
\end{equation}
This is a Hilbert space for the norm given by:
\begin{equation}
\| f \|_\alpha^2 = |f (0)|^2 + (\alpha + 1) \int_\D |f ' (z)|^2 \, (1 - |z|^2)^\alpha \, dA (z) < \infty\, .
\end{equation}
The standard Dirichlet space ${\cal D}$ corresponds to $\alpha = 0$; the Hardy space $H^2$ to $\alpha = 1$ and the standard Bergman space to $\alpha = 2$. 
For more general weights, see \cite{KL}.\par 

We denote by ${\cal D}_{\alpha, \ast}$ the subspace of the $f \in {\cal D}_\alpha$ such that $f (0) = 0$. \par
\smallskip

If $\phi$ is a Schur function, one defines its \emph{weighted Nevanlinna counting function} $N_{\phi, \alpha}$ at $w \in \Omega := \phi (\D)$ as the 
number of pre-images of $w$ with the weight $(1 - |z|)^\alpha$: 
\begin{equation}
N_{\phi, \alpha}  (w) = \sum_{\phi (z) = w} (1 - |z|^2)^\alpha \, .
\end{equation}
For $w \in \D \setminus \phi (\D)$, we set $N_{\phi, \alpha}  (w) = 0$. One has $N_{\phi, 1} = N_\phi$ and $N_{\phi, 0} = n_\phi$. \par 

With this notation, recall the change of variable formula:
\begin{equation} \label{change} 
\int_\D F [\phi (z)] \, | \phi ' (z) |^2 \, (1 - |z|^2)^\alpha \, dA (z) = \int_\Omega F (w) \, N_{\phi, \alpha} (w) \, dA (w) \,.
\end{equation}

Denote by $R_{n, j}$, $n \geq 0$, $0 \leq j \leq 2^n - 1$, the Hastings-Luecking windows:
\begin{displaymath}
R_{n, j} = \Big\{z\in \D\,;\ 1 - 2^{-n} \leq |z| < 1 - 2^{-n - 1}\quad \text{and}\quad 
\frac{2j\pi}{2^n} \leq \arg z < \frac{2(j+1)\pi }{2^n}\,\Big\} \, .
\end{displaymath} 

We can now state. 

\goodbreak 

\begin{theorem} \label{lulu} 
Let $\alpha > - 1$. Let $\phi$ be a Schur function and $p > 0$. Then $C_\phi \in S_p ({\cal D}_{\alpha, \ast})$ if and only if:
\begin{equation} \label{lulu sum} 
\sum_{n = 0}^\infty \sum_{j = 0}^{2^n - 1} \bigg[ 2^{n (\alpha + 2)} \int_{R_{n, j}} N_{\phi,\alpha} (w) \, dA (w) \bigg]^{p/2} < \infty \, .
\end{equation}
If $\phi$ is univalent, \eqref{lulu sum} can be replaced by the purely geometric condition: 
\begin{equation} \label{lulu univalent}
\sum_{n = 0}^\infty \sum_{j = 0}^{2^n - 1} \big[  2^{n (\alpha + 2)} A_\alpha (R_{n, j} \cap \Omega) \big]^{p/2} < \infty \, ,
\end{equation} 
where $A_\alpha$ is the weighted measure $dA_\alpha (w) = (\alpha + 1) \, (1 - |w|^2)^\alpha dA (w)$.  
\end{theorem}

\noindent{\bf Remark.} Of course, every operator in a Schatten class is compact, but we may note that condition~\eqref{lulu sum} implies the compactness 
of $C_\phi$, by \cite{Zor}, Theorem~1 (and \cite{JFA}, Proposition~3.3). \par

\medskip

\noindent{\bf Proof of Theorem~\ref{lulu}.} First, we compute $C_\phi^\ast C_\phi$. Let us fix $f$ and $g$ in the Dirichlet space 
$\mathcal{D}_{\alpha, \ast}$. We have:
\begin{align*}
(\alpha + 1) \int_\D & \big( (C_\phi^\ast C_\phi) (f) \big) ' (z) \, \overline{g ' (z)} \, (1 - |z|^2 )^\alpha \, dA (z) 
=\big\langle f \circ \phi , g \circ \phi \big\rangle_{{\cal D}_{\alpha, \ast}} \\
& = (\alpha + 1) \int_\D (f ' \circ \phi) (z) \overline{(g ' \circ \phi) (z)} \, |\phi ' (z)|^2 \, (1 - |z|^2 )^\alpha \, dA (z) .
\end{align*}
By the change of variable formula, we get:
\begin{displaymath}
\int_\D \big( (C_\phi^\ast C_\phi) (f) \big) ' (z) \overline{g' (z)} \, (1 - |z|^2 )^\alpha dA 
= \int_\D f ' (w) \, \overline{g ' w) \, }N_{\phi,\alpha} (w)\, dA (w) \, ,
\end{displaymath}
which is equivalent to:
\begin{displaymath}
\int_\D \big( (C_\phi^\ast C_\phi) (f) \big) ' (z) \, \overline{G (z)} \, (1 - |z|^2 )^\alpha dA (z) 
=\int_\D f ' (w) \, \overline{G (w)} \, N_{\phi,\alpha} (w) \, dA (w) 
\end{displaymath}
for every function $G$ belonging to the weighted Bergman space ${\mathfrak B}_\alpha^2$.\par

That means that $\big( (C_\phi^\ast C_\phi) (f) \big) ' - f ' . N_{\phi, \alpha} / (1 - |w|^2 )^\alpha$ is orthogonal to the weighted Bergman space 
${\mathfrak B}_\alpha^2$. But $\big( (C_\phi^\ast C_\phi) (f)\big) ' \in {\mathfrak B}_\alpha^2$. Hence 
$\big( (C_\phi^\ast C_\phi) (f) \big) ' $ is the orthogonal projection onto ${\mathfrak B}_\alpha^2$ of the function  
$f ' . N_{\phi, \alpha} / (1 - |w|^2 )^\alpha$. Thus (see \cite{Zhu livre}, \S~6.4.1), we obtain that for every $z \in \D$: 
\begin{align*}
\big( (C_\phi^\ast C_\phi) (f) \big) ' (z) 
& = (\alpha + 1) \int_\D \frac{f ' (w)}{ (1 - \bar{w} z)^{\alpha + 2} } \, \frac{N_{\phi, \alpha} (w)} {(1 - |w|^2 )^\alpha} \, (1 - |w|^2 )^\alpha \, dA (w) 
\\
& = (\alpha + 1) \int_\D \frac{f ' (w)} {(1 -\bar{w} z)^{\alpha + 2} }\, d\mu (w) \\ 
& =  (\alpha + 1) \, T_\mu (f ') (z) \,, 
\end{align*}
where $\mu$ is the positive measure $A$ with weight $N_{\phi, \alpha}$ and $T_\mu$ is the Toeplitz operator on ${\mathfrak B}_\alpha^2$ is introduced 
in \cite{Lue} (let us point out that $\alpha$ in \cite{Lue} corresponds to $- (\alpha + 1)$ in our work). \par

In other words, introducing the map $\Delta (h) = h '$, which is an isometry from ${\cal D}_{\alpha, \ast}$ onto ${\mathfrak B}_\alpha^2$, we have 
$\Delta \circ (C_\phi^\ast C_\phi) = T_\mu \circ \Delta$. We have the following diagram:
\begin{displaymath}
\xymatrix{
{\cal D}_{\alpha, \ast} \ar[d]_{\Delta} \ar[r]^{C_\phi^\ast C_\phi} &  {\cal D}_{\alpha, \ast}  \ar[d]^{\Delta} \\ 
{\mathfrak B}_\alpha^2 \ar[r]^{T_\mu} & {\mathfrak B}_\alpha^2 
}
\end{displaymath}
Hence the approximation numbers of $ T_\mu$ (viewed as an operator  on ${\mathfrak B}_\alpha^2$) and the ones of $C_\phi^\ast C_\phi$ (viewed as 
an operator  on ${\cal D}_{\alpha, \ast}$) are the same. In particular, the membership in the Schatten classes are the same and the final result follows from the 
main theorem in \cite{Lue}: $C_\phi \in S_p ({\cal D}_{\alpha, \ast})$ if and only if $C_\phi^\ast C_\phi\in S_{p/2} ({\cal D}_{\alpha, \ast})$ and that 
holds if and only if:
\begin{displaymath}
\sum_{n = 0}^\infty \sum_{j = 0}^{2^n - 1} \big[ 2^{n (\alpha + 2)} \mu (R_{n, j}) \big]^{p /2} < \infty \, .
\end{displaymath}
Hence $C_\phi \in S_p ({\cal D}_{\alpha, \ast})$ if and only if:
\begin{displaymath}
\sum_{n = 0}^\infty \sum_{j = 0}^{2^n - 1}  \bigg[ 2^{n (\alpha + 2)} \int_{R_{n, j}} N_{\phi, \alpha} (w) \, dA (w) \bigg]^{p / 2} < \infty \, , 
\end{displaymath}
and that ends the proof of Theorem~\ref{lulu}. 
\qed 

\medskip

\noindent{\bf Remark.} In the same way, we can obtain other characterizations for ${\cal D}_{\alpha, \ast}$ by using the ones for ${\mathfrak B}_\alpha^2$ 
given in \cite{LueZhu} and \cite{Zhu}: $C_\phi \in S_p ({\mathfrak B}_\alpha^2)$ if and only if 
$N_{\phi, \alpha + 2} (z)  / \big( \log (1/ |z|) \big)^{\alpha + 2} \in L^{p/2} (\lambda)$, where $d \lambda (z) = (1 - |z|^2)^{- 2} d A (z)$ is the 
M\"obius invariant measure on $\D$, and, when $\phi$ has bounded valence and $ p \geq 2$, if and only if 
$(1- |z|^2) / \big( 1 - |\phi (z)|^2 \big) \in L^{p (\alpha + 2)/ 2} (\lambda)$. Such a result can be found in \cite{Xu}. 
\goodbreak 
 
\subsection {Applications} \label{applis} 

We give several applications of the previous theorem. \par
\begin{corollary} \label{SchattenDirichletalphaimpliqueSchattenDirichletbeta} 
Let $- 1 < \alpha \leq \beta$, $p > 0$, and $\phi$ be a Schur function. Then $C_\phi \in S_p ({\cal D}_{\alpha, \ast})$ implies that 
$C_\phi \in S_p ({\cal D}_{\beta, \ast})$. \par
In particular, $C_\phi \in S_p ({\cal D}_\ast)$ implies that $C_\phi \in S_p (H^2)$, which in turn implies that $C_\phi \in S_p ({\mathfrak B}^2)$.
\end{corollary}

\noindent{\bf Proof.} Assume that $C_\phi \in S_p ({\cal D}_{\alpha, \ast})$. Then 
\begin{displaymath}
\sum_{n = 0}^\infty \sum_{ j = 0}^{2^n - 1} \bigg[ 2^{n (\alpha + 2)} \int_{R_{n, j}} N_{\phi, \alpha} (w) \, dA (w) \bigg]^{p/2} <\infty. 
\end{displaymath}
Since, thanks to Schwarz's lemma, $N_{\phi, \beta} (w) \leq N_{\phi, \alpha} (w) (1 - |w|^2)^{\beta - \alpha}$, we have 
\begin{displaymath}
N_{\phi, \beta} (w) \leq (2.2^{- n} )^{\beta - \alpha} N_{\phi, \alpha} (w) \quad \text{for } w\in R_{n, j}.
\end{displaymath}
It follows that 
\begin{displaymath}
\sum_{n = 0}^\infty \sum_{j = 0}^{2^n - 1} \bigg[ 2^{n (\beta + 2)} \int_{R_{n, j}} N_{\phi, \beta} (w) \, dA (w) \bigg]^{p/2} < \infty \, ,
\end{displaymath}
and that proves Corollary~\ref{SchattenDirichletalphaimpliqueSchattenDirichletbeta}.
\qed 
 
\bigskip

It is known (\cite{JFA}) that composition operators on $H^2$ separate Schatten classes, but the difficulty is that we must not only control the shape of 
$\phi (\partial{\D})$, but also the parametrization $t \mapsto \phi (\e^{it})$, even if $\phi$ is univalent. In the case of the Dirichlet space, this difficulty 
disappears, because only the areas come into play, and we can easily prove the following result. 

\goodbreak 
\begin{theorem} \label{separation Schatten} 
The composition operators on $\mathcal{D}_\ast$ separate Schatten classes, in the following sense. 
Let $0 < p_1 <\infty$. Then, there exists a symbol $\phi$ such that:
\begin{displaymath} 
C_\phi \in \Big( \bigcap_{p > p_1} S_p ({\cal D}_\ast) \Big) \setminus S_{p_1} ({\cal D}_\ast) \, . 
\end{displaymath} 
Similarly, there exists a symbol $\phi$ such that:
\begin{displaymath} 
C_\phi \in S_{p_1} ({\cal D}_\ast) \setminus \Big(\bigcup_{p < p_1} S_p ({\cal D}_\ast) \Big) \, .
\end{displaymath} 
\par
In particular, for every $0 < p_1 < p_2 < \infty$, there exists  $\phi$ such that $C_\phi \in S_{p_2} ({\cal D}_\ast) \setminus S_{p_1} ({\cal D}_\ast)$. 
\end{theorem}

\noindent{\bf Proof.} Let $(h_n)_{n \geq 1}$, with $0 < h_n < 1$, be a sequence of real numbers with limit $0$ to be adjusted, and $J$ the Jordan curve 
formed by the segment $[0, 1]$ and the north and (truncated) north-east sides of the curvilinear rectangles 
\begin{displaymath} 
\{1 - 2^{-n} \leq | z | < 1 - 2^{- n - 1} \} \times \{0 \leq \arg z < 2^{- n} h_n\}. 
\end{displaymath} 
Let $\Omega_0$ be the interior of $J$ and $\Omega = \Omega_0 \cup D (0, 1/8)$. Let $\phi \colon \D \to \Omega$ be a Riemann map such that 
$\phi (0) = 0$. Since $\phi$ is univalent and bounded, it defines a symbol on $\mathcal{D}_\ast$, and the necessary and sufficient 
condition~\eqref{lulu univalent} for membership in $S_{p} ({\cal D}_\ast)$ reads: 
\begin{equation} \label{condition} 
\sum_{n = 0}^\infty [4^n 4^{- n}h_n]^{p/2} = \sum_{n = 0}^\infty {h_n}^{p/2}<\infty. 
\end{equation}
Indeed, it is clear that, for fixed $n$, the Hastings-Luecking windows $R_{n,j}$ satisfy:
\begin{displaymath} 
R_{n, 0} \cap \Omega \neq \emptyset; \quad R_{n,j} \cap \Omega = \emptyset \text{ for } 1\leq j <2^n .
\end{displaymath} 
Therefore, only the Hastings-Luecking windows $R_{n, 0}$ matter. Since:
\begin{displaymath} 
A (R_{n, 0} \cap \Omega) = \iint_{1 - 2^{- n} \leq r < 1 - 2^{- n - 1}, \  0 \leq \theta < 2^{- n} h_n}  r \, dr \, d\theta \approx 4^{- n} h_n \, ,
\end{displaymath} 
we can test the criterion \eqref{condition}. Now, it is enough to take $h_n = (n + 1)^{- 2 /p_1}$ to get: 
\begin{displaymath} 
C_\phi \in \Big( \bigcap_{p > p_1}S_p ({\cal D}_\ast) \Big) \setminus S_{p_1} ({\cal D}_\ast) \, . 
\end{displaymath} 
Similarly, the choice $h_n = (n + 1)^{- 2 /p_1} [\log (n + 2)]^{- 4 / p_1}$, gives a symbol $\phi$ such that: 
\begin{displaymath} 
C_\phi \in S_{p_1} ({\cal D}_\ast) \setminus \Big( \bigcup_{p < p_1} S_p ({\cal D}_\ast) \Big) \, .
\end{displaymath} 
This ends the proof. 
\qed 
 
\medskip \goodbreak

T. Carroll and C. Cowen (\cite{Carroll-Cowen} proved, but only for $\alpha > 0$, that there exist compact composition operators on ${\cal D}_\alpha$ which 
are in no Schatten class (see also \cite{Jones}). In the next result, we shall see that this still true for $\alpha = 0$.

\goodbreak 
\begin{theorem} \label{no Schatten} 
There exists a Schur function $\phi$ such that $C_\phi$ is compact on ${\cal D}_\ast$, but in no Schatten class $S_p ({\cal D}_\ast)$. 
\end{theorem}

\noindent{\bf Proof.} It suffices to use the proof of Theorem~\ref{separation Schatten} and to take, instead of the above $h_n$, $h_n = 1 / \ln (n + 2)$. 
\qed 

\bigskip

For the next application, which will be used in Section~\ref{section capa}, we need to recall the definition of the cusp map $\chi$, introduced in 
\cite{LLQR Bergman}, and later used, with a slightly different definition in \cite{LIQUEFRODR}. Actually, we have to modify it slightly again in order to have 
$\chi (0) = 0$. We first define:
\begin{displaymath}
\chi_0 (z) = \frac{\displaystyle \Big( \frac{z - i}{i z - 1} \Big)^{1/2} - i} {\displaystyle - i \, \Big( \frac{z - i}{i z - 1} \Big)^{1/2} + 1} \, \raise 1pt \hbox{,} 
\end{displaymath}
then: 
\begin{displaymath}
\chi_1 (z) = \log \chi_0 (z), \quad \chi_2 (z) = - \frac{2}{\pi}\, \chi_1 (z) + 1, \quad \chi_3 (z) = \frac{a}{\chi_2 (z)}  \, \raise 1pt \hbox{,} 
\end{displaymath}
and finally:
\begin{displaymath}
\chi (z) = 1 - \chi_3 (z) \, ,
\end{displaymath}
where $a = 1 - \frac{2}{\pi} \log (\sqrt{2} - 1) \in (1, 2)$ is chosen in order that $\chi (0) = 0$. The image $\Omega$ of the (univalent) cusp 
map is formed by the intersection of the  inside of the disk $D \big(\frac{a}{2} \raise 1pt \hbox{,} \frac{a}{2} \big)$ and the outside of the two disks  
$D \big(\frac{i a}{2} \raise 1pt \hbox{,} \frac{a}{2} \big)$ and $D \big( - \frac{ i a}{2} \raise 1pt \hbox{,} \frac{a}{2} \big)$. 

\goodbreak 
\begin{corollary} \label{pointe} 
If $\chi$ is the cusp map, then $C_\chi$ belongs to all Schatten classes $S_p ({\cal D}_\ast)$, $p > 0$.
\end{corollary} 

\noindent{\bf Proof.} Since $\chi$ is univalent, $\chi (0) = 0$, and $\Omega = \chi (\D)$ has finite area, we have $\chi \in {\cal D}_\ast$. 
A little elementary geometry shows that, for some constant $C$, we have:
\begin{equation} \label{geoana} 
w \in \Omega, \  0 < h < 1 \text{ and } | w | \geq 1 - h \quad \Longrightarrow \quad | \Im w | \leq C h^2.
\end{equation}
It follows (changing $C$ if necessary) that $R_{n, j} \cap \Omega$ is contained in a rectangle of sizes $2^{- n}$ and $C \, 4^{- n}$ and with area $C \, 8^{-n}$. 
Hence, for a given $n$, at most $C$ of the Hastings-Luecking windows $R_{n, j}$ can intersect $\Omega$. Therefore, the series in Theorem~\ref{lulu} reduces, 
up to constants, to the series:
\begin{displaymath}
\sum_{n = 0}^\infty (4^n 8^{- n} )^{p/2} = \sum_{n = 0}^\infty 2^{- n p} \, ,
\end{displaymath}
which converges for every $p > 0$. 
\qed 

\goodbreak 


\section{Logarithmic capacity and set of contact points} \label{section capa}  

In view of the result of  \cite{Gallardo} mentioned in the introduction, if ${ \rm Cap} \, K > 0$, there is no hope to find a symbol $\phi$ such that 
$E_\phi = K$ and $C_\phi$ is Hilbert-Schmidt on $\mathcal{D}_\ast$. But as was later proved in \cite{ELFKELSHAYOU}, ${ \rm Cap} \, K > 0$ is the only 
obstruction. We can improve on the results from \cite{ELFKELSHAYOU} as follows: our composition operator is not only Hilbert-Schmidt, but in any 
Schatten class; moreover, we can replace $E_\phi = K$ by $E_\phi = E_\phi (1) = K$. 
\goodbreak

\begin{theorem} \label{main th} 
For every compact set $K$ of the unit circle $\T$ with  logarithmic capacity ${\rm Cap} \, K = 0$, there exists a Schur function $\phi$ with the following 
properties: \par
\smallskip

1) $\phi \in A (\D) \cap {\cal D}_\ast := A$, the ``Dirichlet algebra''; \par
\smallskip

2) $E_\phi = E_\phi (1) = K$; \par
\smallskip

3) $C_\phi \in  \bigcap_{p > 0} S_p ({\cal D}_\ast)$. \par 
\smallskip

\noindent In fact, the approximation numbers of $C_\phi$ satisfy $a_n (C_\phi) \leq a \exp(- b \sqrt {n})$.
\end{theorem}

This theorem actually results of the particular following case and the properties of the cusp map seen in Section~\ref{applis}.  

\begin{theorem} \label{q peak}
For every compact set $K \subseteq \partial{\D}$ of logarithmic capacity ${\rm Cap} \, K = 0$, there exists a Schur function $q\in A (\D) \cap {\cal D}_\ast$ 
which peaks on $K$ and such that the composition operator $C_q \colon {\cal D}_\ast \to {\cal D}_\ast$ is bounded (and even Hilbert-Schmidt).
\end{theorem}

Recall that a function $q \in A (\D)$, the disk algebra, is said to \emph{peak} on a compact subset $K \subseteq \partial{\D}$ (and is called a 
\emph{peaking function} ) if: 
\begin{displaymath}
q (z) = 1 \text{\ if } z\in K \,; \quad | q(z) | < 1 \text{ if } z \in \overline{\D} \setminus K \, .
\end{displaymath}

\noindent{\bf Proof of Theorem~\ref{main th}.} We simply take for $\phi$ the composed map $\phi = \chi \circ q$, where $\chi $ is the cusp map and $q$ 
our peaking function. Recall that $\chi \in A (\D)$ and that $\chi$ peaks on $\{1\}$. We take advantage of this fact by composing with $q$, for which 
$C_q \colon {\cal D}_\ast \to {\cal D}_\ast$ is bounded as well as $C_\chi$ (since $\chi$ is univalent). We clearly have $\phi \in A (\D)$, 
$\phi (z) = \chi (1) = 1$ for $z \in K$, and $|\phi (z) | < 1$ for $z\notin K$, since then $| q (z) | < 1$. Therefore $E_\phi (1) = K$. Moreover,  
$C_\phi$ being bounded on ${\cal D}_\ast$, we have in particular $\phi = C_\phi (z) \in {\cal D}_\ast$. Since $C_\phi = C_q\circ C_\chi$, we get 3), by   
 Corollary~\ref{pointe}. \par
In \cite{approx Dirichlet}, we prove that $a_n (C_\chi) \leq a \exp(- b \sqrt {n})$. Since $a_n (C_\phi) \leq \| C_q \| \,  a_n (C_\chi)$, by the ideal property 
of approximation numbers, this ends the proof of Theorem~\ref{main th}.  
\qed
\bigskip

In turn, the proof of Theorem~\ref{q peak} relies on the following crucial lemma. 
\goodbreak

\begin{lemma} \label{crucial} 
Let $K \subseteq \partial{\D}$ be a compact set such that  ${\rm Cap} \, K = 0$. Then, there exists a function 
$U \colon \overline{\D} \to \R^{+} \cup \{\infty\}$, such that: \par
\smallskip

1) $U (z) = \infty$ if and only if $z \in K$; 
\par

2)  $U \geq 1$ on $\overline{\D}$; \par
\smallskip

3) $U$ is continuous on $\overline{\D} \setminus K$, harmonic in $\D$ and $\int_{\D} | \nabla U |^2 \, dA < \infty$;  \par 
\smallskip 
 
4) $\lim_{z \to K,\ z\in \overline{\D}} U (z) = \infty$; \par
\smallskip
 
5) the conjugate function $V = \tilde U$ is continuous on $\overline{\D} \setminus K$.
\end{lemma}

\noindent{\bf Proof of Theorem~\ref{q peak}.} Taking this lemma for granted, let us end the proof of the theorem. We set $f = U +i V$, $a = \e^{- 1 / f (0)}$ 
and $q = \phi_a \circ \e^{- 1/ f} $, where $\phi_a (z) = \frac{z - a}{1 - \bar{a} z}$. In view of the third and fourth items of the lemma, we have $q \in A (\D)$. 
Since $U \geq 1$, Lemma~\ref{CS HS} shows that $C_q$ is Hilbert-Schmidt on ${\cal D}_\ast$. Moreover, for $z \in K$, one has $f (z) = \infty$ and hence 
$q (z) =1$ since $\phi_a (1) = 1$ because $a \in \R$ (since $f (0) = U (0)$). On the other hand, when  $z \notin K$, one has $| f( z) | < \infty$ and hence 
$| q (z) | < 1$. Therefore $q$ peaks on $K$. 
\qed 

\bigskip

\noindent{\bf Proof of Lemma~\ref{crucial}.} This proof is strongly influenced by that of Theorem~III, page~47, in \cite{KAHSAL}. Let: 
\begin{equation}
L (z) = \log \Big(\frac{\e}{1 - z} \Big) = P (z) + i \, Q (z),
\end{equation}
with 
\begin{displaymath}
\quad P (z) = \log \frac{\e}{| 1 - z |} \text{ and } Q (z) = - \arg (1 - z), \ | Q (z) | \leq \frac{\pi}{2}, \quad  z \in \overline{\D} \setminus \{1\} \, ,
\end{displaymath}
and write: 
\begin{displaymath}
P (z) \sim \sum_{n \in \Z} \gamma_n \, z^n \,,
\end{displaymath}
with 
\begin{displaymath}
\gamma_n = 1 / (2 \, |n|) \quad  \text{if  } n \neq 0 \,,  \quad \text{and }  \gamma_0 = 1  \, .
\end{displaymath}

For $0 < \eps < 1 / 2$, let $K_\eps =\{ z \in \T \, ; \  {\rm dist} \, (z, K) \leq \eps \}$, $\mu_\eps$ its equilibrium measure, and $U_\eps$ the logarithmic 
potential of $\mu_\eps$, that is:
\begin{displaymath}
U_\eps (z) = \int_{K_\eps} \log \frac{\e}{|z - w|} \, d\mu_\eps (w)  \, , 
\end{displaymath}
that we could as well write (since $K_\eps \subseteq \T$):
\begin{displaymath} 
U_\eps (z) = \int_{K_\eps} P (z \, \bar{w}) \, d\mu_\eps (w) \, .
\end{displaymath} 
\smallskip 

Let us set: 
\begin{equation}  
f_\eps (z) = \int_{K_\eps} L (z \, \bar{w} ) \, d\mu_\eps (w) = U_\eps (z) + i V_\eps (z) \,,
\end{equation} 
with
\begin{displaymath}
V_\eps (z) = \int_{K_\eps} Q (z \, \bar{w} ) \, d\mu_\eps (w) \, .
\end{displaymath}
Then, if $I_\eps$ is the energy of $\mu_\eps$, one has (see \cite{Ross}, Section~4) 
$I_\eps = 1 + \sum_{n = 1}^\infty \frac{|\widehat{\mu_\eps} (n)|^2}{n}$, where 
$\widehat{\mu_\eps} (n) = \int_\T {\overline w \,}^n \, d\mu_\eps (w)$ is the $n$-th Fourier coefficient of $\mu_\eps$, and: 
\begin{equation} \label{observe} 
f_\eps \in {\cal D} \quad \text{and} \quad \| f _\eps \|_{\cal D}^2 = I_\eps \, .
\end{equation} 
Note that $\| f _\eps \|_{\cal D} \geq 1$. \par
\smallskip 

We claim that there exist $\delta > 0$ and $0 < r < 1$ such that: 
\begin{equation} \label{unif}
z \in \overline{\D} \text{ and } {\rm dist} \, ( z, K) \leq \delta \quad \Longrightarrow \quad U_\eps (r z) \geq I_\eps / 2 
\end{equation} 

Indeed, let $P_a (t) = \frac{1 - |a|^2}{|\e^{it} - a|^2}$ be the Poisson kernel at $a \in \D$. Since $U_\eps$ is harmonic in $\D$ and integrable on $\T$ 
(\cite{Conway}, Proposition~19.5.2), one has, for every $z \in \D$: 
\begin{equation} \label{Poisson}
U_\eps (z) = \int_{- \pi}^{\pi} U_\eps ( \e^{it})  \, P_z (t) \, \frac{dt}{2 \pi} \, .
\end{equation} 

Let now $\delta \leq \eps / 4$, to be adjusted later, and take $1 - \delta \leq r < 1$. Suppose that ${\rm dist} (z, K) \leq \delta$, with $z \in \overline{\D}$, 
and let $u \in K$ such that $|z - u | \leq \eps / 4$. Note that then $|r z - u | \leq (1 - r) + |z - u| \leq \eps / 2$. It follows from \eqref{Poisson} that: 
\begin{displaymath} 
I_\eps - U_\eps (r \, z) = \int_{- \pi}^{\pi} [ I_\eps - U_\eps ( \e^{it})] \, P_{r z} (t) \, \frac{dt}{2 \pi} 
\end{displaymath} 
(it is useful to recall that $U_\eps (z) \leq I_\eps$ for every $z \in \C$). Set:
\begin{displaymath} 
J_1 = \int_{|\e^{it} - r z| \leq \eps / 2} [ I_\eps - U_\eps ( \e^{it})] \, P_{r z} (t) \, \frac{dt}{2 \pi} 
\end{displaymath} 
and
\begin{displaymath} 
J_2 = \int_{|\e^{it} - r z| > \eps / 2} [ I_\eps - U_\eps ( \e^{it})] \, P_{r z} (t) \, \frac{dt}{2 \pi} \, .
\end{displaymath} 

For the integral $J_1$, we have:
\begin{displaymath} 
|\e^{it} - u | \leq |\e^{it} - r z  | + |r z - u | \leq \eps \, ; 
\end{displaymath} 
therefore $\e^{it} \in K_\eps$. Since $U_\eps = I_\eps$ Lebesgue-almost everywhere on $K_\eps$, by Frost\-man's Theorem, we get $J_1 = 0$. \par

For the integral $J_2$, we have:
\begin{displaymath} 
P_{r z} (t)  \leq \frac{2 (1 - r \, |z|)}{(\eps / 2)^2} \leq 2\, \frac{(1 - r) + r (1 - |z|)}{(\eps / 2)^2} \leq \frac{4 \delta}{(\eps / 2)^2} 
= \frac{16 \delta}{\eps^2} \, ;
\end{displaymath} 
hence (since $U_\eps (\e^{it}) \geq 0$):
\begin{displaymath} 
J_2 \leq \frac{16 \delta}{\eps^2}\, I_\eps \, .
\end{displaymath} 

Therefore, if we choose $0 < \delta \leq \eps^2 / 32$, we get:
\begin{displaymath} 
0 \leq I_\eps - U_\eps (r \, z) \leq I_\eps / 2 \,, 
\end{displaymath} 
which gives \eqref{unif}. \qed
\par
\medskip

Now, as ${\rm Cap} \, K = 0$, we know from \eqref{explose} that $\lim_{\eps \to 0^+} I_\eps = \infty$, and we can adjust a sequence $\eps_j \to 0^+$  so 
that:
\begin{equation} \label{adjust} 
I_{\eps_j}  \geq 4 \, j^6 \, .
\end{equation}

Using \eqref{unif}, we find two sequences $(\delta_j)_j$ and $(r_j)_j$, with $0 < \delta_j \to 0$ and $1 > r_j \to 1$, such that, for every $j \geq 1$, 
\begin{equation} \label{adjustencore}  
z \in \overline{\D} \text{ and } {\rm dist} \, (z, K) \leq \delta_j \quad \Longrightarrow \quad U_{\eps_j} (r_j z) \geq I_{\eps_j}/2  .
\end{equation}

Finally, let us set:
\begin{equation}
f_j (z) = f_{\eps_j} (r_j z) 
\end{equation}
and 
\begin{equation} \label{choice} 
f = U + i V = 1+ \sum_{j = 1}^\infty j^{- 2} \frac{f_j \ }{\| f_j  \|_{\cal D}} \, \cdot
\end{equation}
The series defining $f$ is absolutely convergent in $\mathcal{D}$. Note that $f (0)$ is real. \par 
\smallskip

We now have: \par
\smallskip

1) $f$ is continuous on $\overline{\D} \setminus K$. \par

Indeed, let $z \in \overline{\D} \setminus K$. Then, ${\rm dist} \, (z, K) > 0$ and there exists a neighbourhood $\omega$ of $z$ in $\overline{\D}$, an 
integer $j_0 = j_0 (z)$ and a positive number $\delta > 0$ such that:
\begin{displaymath}
w \in \omega \text{ and }  j \geq j_0  \quad \Longrightarrow \quad {\rm dist} \, (r_j w, K_{\eps_j}) \geq \delta. 
\end{displaymath}
We then have, for $w \in \omega$ and $j \geq j_0$: 
\begin{align*}
| f_{\eps_j} (w) | 
& =\bigg| \int_{K_{\eps_j}} \log \frac{\e}{r_j w - u} \, d\mu_{\eps_j} (u) \bigg| \\ 
& \leq \int_{K_{\eps_j}} \Big( \log \frac{\e}{| r_j w - u |} + \frac{\pi}{2} \Big) \, d\mu_{\eps_j} (u) 
\leq \log \frac{\e}{\delta} +\frac{\pi}{2} := C \, ,
\end{align*}
since $\mu_{\eps_j}$ is a probability measure supported by $K_{\eps_j}$. Therefore, the series defining $f$ is normally convergent on $\omega$ since its general term is dominated by $j^{- 2} C $ on $\omega$. Since the functions $f_j$ are continuous on $\overline{\D}$, this shows that $f$ is continuous at $z$. 
\par
\smallskip

2) $U (z) := \Re f (z) \geq 1$. \par 
This is obvious since, for every $z\in \overline{\D}$, 
\begin{displaymath}
U_\eps (z) := \Re f_\eps (z) = \int_{K_\eps}  \log \frac{\e}{|z - u |} \, d\mu_\eps (u) \geq 0 \, . 
\end{displaymath}

3) $\lim_{z \to K, z \in \overline{\D}} U (z) =\infty$. \par \smallskip

Indeed,  let $A > 0$. Take an integer $j \geq A$ and suppose that ${\rm dist} \, (z, K) \leq \delta_j$. Then, using the positivity of the $U_{\eps_k}$'s as well 
as \eqref{observe}, \eqref{adjust} and \eqref{adjustencore}, we have:
\begin{displaymath}
U (z) \geq j^{- 2} \frac{U_{\eps_j} (r_j z)} {\| f_{\eps_j} \|_{\cal D}} 
\geq j^{- 2} \frac{I_{\eps_j}/2 } {\sqrt{I_{\eps_j} } } \geq j \geq A. 
\end{displaymath}
This ends the  proof of our claims, and of Lemma~\ref{crucial}. 
\qed

\bigskip

To end this paper, let us mention the following version of the classical Rudin-Carleson Theorem. Though it is not the main subject of this paper, it has the 
same flavor as Theorem~\ref{q peak}. We do not give a proof, but only mention that it can be obtained by mixing the proofs of Theorems~III.E.2 and III.E.6 
in \cite{Woj} (see pages 181--187). 

\begin{theorem} 
Let $K$ be a compact subset of $\T$ with ${\rm Cap}\, K = 0$. Given any continuous strictly positive function $s \in C (\T)$ equal to $1$ on $K$, we can 
find, for every $h \in C (K)$  and every $\eps > 0$, a function $f \in A (\D) \cap {\cal D}$ such that $f_{\mid K} = h $ and:
\begin{displaymath} 
\qquad  |f (\theta) | \leq (1 + \eps) \, \| h \|_\infty\, s (\theta) \, ,\ \forall \theta \in \T \ ; 
\qquad \| f \|_{\cal D} \leq (1 + \eps) \, \| h \|_\infty \, . 
\end{displaymath} 
\end{theorem} 
%


\bigskip

\vbox{\small \noindent{\it
{\rm Pascal Lef\`evre}, Univ Lille Nord de France, U-Artois, \\
Laboratoire de Math\'ematiques de Lens EA~2462 {\rm \&} F\'ed\'eration CNRS Nord-Pas-de-Calais FR~2956, \\
Facult\'e des Sciences Jean Perrin, Rue Jean Souvraz, S.P.\kern 1mm 18, \\ 
F-62\kern 1mm 300 LENS, FRANCE \\ 
pascal.lefevre@euler.univ-artois.fr
\smallskip

\noindent
{\rm Daniel Li}, Univ Lille Nord de France, U-Artois, \\
Laboratoire de Math\'ematiques de Lens EA~2462 {\rm \&} F\'ed\'eration CNRS Nord-Pas-de-Calais FR~2956, \\
Facult\'e des Sciences Jean Perrin, Rue Jean Souvraz, S.P.\kern 1mm 18, \\ 
F-62\kern 1mm 300 LENS, FRANCE \\ 
daniel.li@euler.univ-artois.fr
\smallskip

\noindent
{\rm Herv\'e Queff\'elec}, Univ Lille Nord de France, \\
USTL, Laboratoire Paul Painlev\'e U.M.R. CNRS 8524, {\rm \&}  F\'ed\'eration CNRS Nord-Pas-de-Calais FR~2956, \\
F-59\kern 1mm 655 VILLENEUVE D'ASCQ Cedex, 
FRANCE \\ 
Herve.Queffelec@univ-lille1.fr
\smallskip

\noindent
{\rm Luis Rodr{\'\i}guez-Piazza}, Universidad de Sevilla, \\
Facultad de Matem\'aticas, Departamento de An\'alisis Matem\'atico {\rm \&} IMUS,\\ 
Apartado de Correos 1160,\\
41\kern 1mm 080 SEVILLA, SPAIN \\ 
piazza@us.es\par}
}

\end{document}